\documentclass{article}[11 pts]
\usepackage{url}
  \usepackage{graphicx}
  \usepackage{psfrag}
  \usepackage{latexsym}
  \usepackage{amsfonts}
\newtheorem{thm}{Theorem}[section]
\newtheorem{prop}[thm]{Proposition}

\newtheorem{lemma}[thm]{Lemma}

\newcommand{\Z}{\mathbb{Z}}
\newcommand{\R}{\mathbb{R}}
\newcommand{\Q}{\mathbb{Q}}
\newcommand{\N}{\mathbb{N}}
\newcommand{\proof}{\noindent {\it Proof: }}
\newcommand{\qed}{\vrule height4pt width4pt depth4pt}
\begin{document}
\bibliographystyle{alpha}
\unitlength=1cm
\title{\bf Algebraic Unimodular Counting}
\author{
    \begin{tabular}[t]{c}
      Jes\'us A. De Loera
      \\[-.1cm] {\footnotesize Dept. of Mathematics}
      \\[-.1cm] {\footnotesize Univ. of California, Davis}
      \\[-.1cm] {\footnotesize deloera@math.ucdavis.edu}
    \end{tabular}
    \begin{tabular}[t]{c}
      Bernd Sturmfels
      \\[-.1cm] {\footnotesize Dept. of Mathematics}
      \\[-.1cm] {\footnotesize Univ. of California, Berkeley}
      \\[-.1cm] {\footnotesize bernd@math.berkeley.edu}
    \end{tabular}
}
\date{}

\maketitle

\begin{abstract}
\noindent   We study algebraic algorithms for expressing the number
   of non-negative integer solutions to a unimodular system
   of linear equations as a function of the right hand side.
   Our methods include Todd classes of toric varieties via
   Gr\"obner bases, and rational generating functions as in
   Barvinok's algorithm. We report polyhedral and computational results for
   two special cases:  counting contingency tables and
   Kostant's partition function.
\end{abstract}

\section{Introduction}
The object of study in this paper is  the vector partition function
 $$ \phi_A(b) \quad = \quad
\# \, \bigl\{\,x \,:\, Ax=b, x \geq 0 , \,\, x \, \,{\rm integral} \,
\bigr\} , $$
where $A$ is a fixed $d \times n$-unimodular integer matrix and $b$
is a variable vector in $\Z^d$. Here we say that $A$ is
{\em unimodular} if the polyhedron $\{x : Ax=b, x \geq 0\}$ has only
integral vertices whenever $b$ is in the lattice spanned by the
columns of $A$. This is a slight generalization of
the definition of ``unimodular'' used in
\cite[\S 19]{schrijver}. We further assume that
$\,Ker(A) \cap \R_{\geq 0}^n = {0} $, which is equivalent to
$\,\phi_A(b) < \infty\,$ for all $b$.
We regard $\phi_A$ as a function on $\,{\rm cone}(A)$, the cone
of non-negative linear combinations of   the columns of $A$, since
$\phi_A(b) = 0$ if $b$ is not in  $\, {\rm cone}(A)$.
The following result about  vector partition functions
is well-known (see e.g.~\cite{sturmfels}):

\begin{thm} \label{structure}
The function $\phi_A$ is piecewise polynomial of degree $n-rank(A)$.
Its domains of polynomiality are convex polyhedral cones, called
{\em chambers} of $A$.
\end{thm}

The main purpose of this paper is to develop practical methods for
unimodular counting. By {\it unimodular counting} we mean
preprocessing the given unimodular matrix $A$ and generating the
polynomials for $\phi_A$ on the various chambers. Each output
polynomial is represented either explicitly as a sum of monomials, or
implicitly as an oracle which allows for quick evaluation of $\phi_A$
at any $b$ in that chamber.  Unimodular counting has many
applications, ranging from statistics \cite{dinwoodie} and randomized
algorithms \cite{welsh} to representation theory \cite{kirillov,
schmidtbincer}. For instance, the widely known problem of counting
{\it contingency tables} is the case when $A$ is the incidence
matrix of a complete bipartite graph \cite{persi2, persibernd}.

Our benchmark on unimodular counting is the work of Mount \cite{mount,
mount2}.  His approach is based on interpolating the chamber
polynomials, by evaluating $\phi_A(b)$ for sufficiently many right
hand sides $b$, coupled with divide-and-conquer decompositions and
advanced parallel computation techniques. Both the evaluation and the
divide-and-conquer schemes depend on the specific matrix $A$.  Mount
reports the complete solution for contingency tables of size $4 \times
4$. In Welsh's survey \cite{welsh} on approximate counting, Mount's
computations for $4 \times 4$-tables are mentioned as the limit for
exact counting on today's computers.

Mount's method does not take full advantage of the rich algebraic
structure underlying $\phi_A$. On page 64 of his thesis \cite{mount},
he writes {\sl ``There are some results in commutative algebra that
relate the (chamber) polynomials to ``Hilbert series'' and ``Todd
classes'', but these structures encode a lot of information and are in
themselves hard to compute.  The strategy taken here is to assume
access to a counting oracle .... and then recover the desired
polynomial by interpolating...''}

We shall demonstrate that algebraic algorithms perform much
better than Mount had surmised. In fact, using rather simple test
implementations, we can now count $4 \times 5$ and $5 \times 5$
contingency tables with arbitrarily large margins.

The algebraic methods apply to any unimodular matrix $A$ and work
independently of the size of the right-hand-side vector $b$. In fact,
our original motivation for this project was the open problem, stated
by Kirillov \cite[page 57]{kirillov}, of computing the number of
chambers for {\sl Kostant's partition function} of the root system
$A_{m-1}$. In this special case, our unimodular matrix is the
incidence matrix of the complete directed graph $K_m$.  We solve
Kirillov's problem for $m \leq 7$, and we compute all chamber
polynomials up to $m \leq 6$.  Using these polynomials we provide an
on-line calculator for Kostant's partition function at
\url{www.math.ucdavis.edu/~deloera/kostant.html}.  We also
prove some other new results on the geometry of chamber complexes of
unimodular matrices.

This paper presents two algebraic algorithms for unimodular counting:
\begin{enumerate}
\item A {\it Gr\"obner bases algorithm}, which computes the Todd class
of the toric variety defined by our polytope $\,\{x \geq 0 : A x = b\}$,
is given in Section 2.
\item The {\it BBKLP method}, which computes the generating function for
all lattice points in the polytope $\{x \geq 0 : A x = b\}$,
is given in Section 3.
\end{enumerate}
The acronym BBKLP refers to five mathematicians: Barvinok, Brion,
Khovanskii, Lawrence, and Pukhlikov. The most important complexity
result in our area is Barvinok's polynomial-time algorithm for
counting lattice points in rational polytopes of fixed dimension
\cite{barvinok, dyerkannan}. Barvinok's algorithm is based on earlier
work by Brion, Khovanskii, Lawrence, and Pukhlikov. For a complete
bibliography see the survey article of Barvinok and Pommersheim
\cite{BarviPom}.  When $A$ is unimodular, Barvinok's algorithm
specializes to the BBKLP method and runs very fast in practice.
This answers a question of Mount \cite[page 56]{mount}.

We implemented methods (1) and (2) in the computer algebra packages
{\tt Macaulay 2} and {\tt Maple} respectively.  Details are described
in Sections 4 and 5. We expect a significant further speed-up by
combining our algebraic approach with Mount's parallel computing
techniques. In a future project we will extend the various methods for
computing $\phi_A$ to non-unimodular matrices $A$.

\section{Method One: Counting Using Gr\"obner bases}

We describe now our first algebraic algorithm for solving the following
counting problem associated with any unimodular $d \times n$-matrix
$A$: {\sl Determine the number $\phi_A(b)$ of non-negative integer
solutions $u \in \N^n$ of the linear equations $\, A \cdot u =
b$.}

The following discussion makes use of well-known facts from algebraic
geometry (see \cite{Ful}); specifically, we demonstrate how to
effectively compute the Todd cohomology class of a toric manifold
defined by a unimodular matrix.

Our running example is the following unimodular
$3 \times 5$-matrix:
$$ A \quad = \quad \pmatrix{ 1 & 0 & 0 & 1 & 1 \cr 0 & 1 & 0 & 1 & 0
\cr 0 & 0 & 1 & 0 & 1 \cr}. $$
The vector partition function for this matrix equals {\small
$$ \phi_A(a,b,c) \, = \,
\cases{
b c+b+c+1 & if $a \geq b+c $ and $b,c \geq 0$, \cr
{1 \over 2} a^2+ {3 \over 2} a+1 & if $\,
{\rm min}\{b,c \} \geq a \geq 0 $, \cr
ab -{1 \over 2} b^2 + {1 \over 2} b+a+1
& if $\,c \geq a \geq b \geq 0$,\cr
ac -{1 \over 2} c^2 + {1 \over 2} c+a+1 &
if $\,b \geq a \geq c $,  \cr
ab \! + \!  ac \!- \! {1 \over 2}(a^2 \! + \! b^2 \! + \!  c^2)
\!  + \! {1 \over 2} (a \! + \! b \! + \! c) \! + \! 1 &
if $ b+c \geq a \geq {\rm max} \{b,c\} $.  \cr}
$$}
For our exposition it is more
convenient to express the vector partition function as
$\,\psi_A  : \N^n \rightarrow \N\,$
where $\psi_A(v)$ is the number of solutions
$u \in \N^n$ to the equation $A u = A v$. Clearly,
$\psi_A$ and $\phi_A$ are related by  a simple
transformation. For instance, in our example we have
$\psi_A(a,b,c,d,e) \, = \, \phi_A(a+d+e, b+d, c+e) $.

The {\it chamber complex} of a unimodular matrix $A$ is defined
as the common refinement of all  triangulations of $A$.
For the $3 \times 5$-matrix $A$ above, the chamber complex is
the given subdivision of $\,{\rm cone}(A) = \R_{\geq 0}^3 \,$ into five
triangular cones.
We refer to \cite{billeraetal} and
\cite{deloeraetal} for details on chamber complexes and \cite{ziegler}
for an introduction to Gale transforms and triangulations.
We assume that $\,{\rm rank}(A) = d$.

\begin{lemma} The chambers of $A$
are in bijection with the regular triangulations of
any  Gale transform $\hat{A}$ of $A$.
Non-regular triangulations of $\hat{A}$ are in
bijection with the virtual chambers of $A$.
\end{lemma}

Thus generating the chambers of our unimodular matrix $A$ is the same
as generating all regular triangulations of a Gale transform
$\hat{A}$. It is well-known that the regular triangulations can be
generated by a applying bistellar flips to a seed regular
triangulation (see \cite{ziegler}). Bistellar flips are topological
operations that transform a triangulation into another. One has to be
careful as sometimes a flip creates non-regular triangulations, but
regularity of a triangulations can be checked by linear programming.
When necessary we have performed these calculations using the software
packages {\tt Puntos} \cite{puntos}  and {\tt Topcom} \cite{rambau}.

We first characterize the chamber complex in algebraic terms.  Let $S
= k[x_1,\ldots,x_n]$ be the polynomial ring over a field $k$ which
contains the rational numbers.  The variables of $S$ index the
columns of the matrix $A = (a_{ij})$.  Let $J_A$ denote the ideal in
$S$ generated by the binomials $\, x_1^{a_{i1}} x_2^{a_{i2}} \cdots
x_n^{a_{in}} \, - \, 1\,$ for $i = 1,2,\ldots,d$. For any positive
weight vector $w \in \R^n$, let $in_w(J_A)$ denote the ideal
generated by the $w$-initial forms of the binomials in $J_A$. If $w$ is
generic, then $in_w(J_A)$ is a monomial ideal. It was shown in
\cite[Corollary 8.9]{Stu}  that the matrix $A$ is unimodular if and
only if all initial monomial ideals $in_w(J_A)$ are square-free.  Two
weight vectors $w$ and $w'$ in $\R^n$ lie in the same cone of
the {\em Gr\"obner fan} if $in_w(J_A) = in_{w'}(J_A)$.  By the
results in  \cite[\S 8]{Stu} this happens  if and
only if, for every linearly independent subset $\sigma = \{a_{i_1},
\ldots,a_{i_r}\}$ of column vectors of $A$, the vector $Aw$ lies in
the cone spanned by $\sigma$ if and only if the vector $Aw'$ lies in
the cone spanned by $\sigma$. This implies the following result:

\begin{prop}
The chamber complex of $A$ equals the
Gr\"obner fan of $J_A$.
\end{prop}

Algebraic algorithms for computing Gr\"obner fans are described in
\cite[\S 3]{Stu}. The state of the art on this subject
is the  work of Huber and Thomas \cite{birkrekha}.
We now explain how to compute the polynomial
representing $\psi_A$ on any given chamber. Suppose
that $w$ is a positive integer vector in the interior
of that chamber. Then $\, M = in_w(J_A) \,$ is
a square-free monomial ideal.
It was shown in \cite[Corollary 7.4]{SWZ}
that $M$ encodes the face poset of the simple polytope
$$\, P_w \quad = \quad \bigl\{ \,u \in \R^n \,\, : \,\, u \geq 0
\,\,\hbox{and} \,\, A u = A w \, \bigr\}. $$
For any
$(n-d)$-element subset $I$ of $\{1,\ldots,n\}$, the equations $\, u_i
= 0, \, i \in I \,$ define a vertex of $P_w$ if and only if $\langle
x_j  :  j \not\in I \rangle $ is a minimal prime of $M$. Writing
$\Sigma_w$ for the normal fan of the simple
polytope $P_w$, this can be restated as follows:

\begin{prop}
The  Stanley-Reisner ideal of the  fan $\Sigma_w$
equals $M = in_w(J_A)$.
\end{prop}

In our running example, with $w=(1,1,1,1,1)$, the polytope $P_w$ is a
pentagon and the fan $\Sigma_w$ has five rays in the plane. This is
encoded by the ideal
\begin{equation}
\label{(3)}
 M \,\, = \,\, \langle A,B,C \rangle \, \cap \, \langle A,B,E
\rangle \, \cap \, \langle B,D,E \rangle \, \cap \, \langle C,D,E
\rangle \, \cap \, \langle A,C,D \rangle .
\end{equation}
Returning to
the general case, our goal is to count the lattice points in the
polytope $P_w$.  We use known methods from toric geometry for this
computation.  An introduction can be found in Section 5.3 in Fulton's
book \cite{Ful}. See also \cite{BarviPom}.

Let $X_w$ denote the projective toric variety defined by the fan
$\Sigma_w$. The variety $X_w$ is smooth, for all generic $w$, since
$A$ is unimodular.  Let $L_A$ denote the ideal in $S =
k[x_1,\ldots,x_n]$ generated by the linear forms $\, b_1 x_1 + \cdots
+ b_n x_n $, where $(b_1,\ldots,b_n)$ runs over all vectors in the
kernel of the matrix $A$.  The cohomology ring of $X_w$ with
coefficients in our field $k$ is the artinian graded $k$-algebra
\begin{equation}
\label{(4)}
  H^*(X_w; k) \quad = \quad \bigoplus_{r=0}^{n-d} H^{2r}(X_w, k)
\quad = \quad S/(M + L_A) .
\end{equation}
Arithmetic operations in this algebra are performed
using normal form reduction relative to any Gr\"obner basis of the
ideal $M+L_A$. Since $X_w$ is an irreducible complex manifold of dimension
$n-d$,  the  top cohomology group $H^{2n-2d}(X_w, k)$ is  a
one-dimensional vector space. There is a canonical choice of
a basis vector  for that one-dimensional $k$-vector space,
namely any
square-free monomial $ \, \prod_{i \in I}  x_i \,$ which indexes a vertex
of $P_w$. This is equivalent to
  $\langle x_j \, : \, j \not\in I \rangle \,$ being a
minimal prime of $M$. Since $X_w$ is smooth, any two such monomials
are congruent to each other modulo $\, M + L_A $. The resulting
 element of  $\, H^*(X_w; k) \,$ represents the
cohomology class which is Poincar\'e dual to a point.

The following rule uniquely defines a $k$-linear functional
called the  {\it integral}:
$$\,  H^*(X_w; k)  \rightarrow k, \,
\, p \,\mapsto \, \int_{X_w} \! p. $$
Writing ${\rm top}(p)$ for the degree $n-d$ component of $p$, we
require that $\, {\rm top}(p)
- (\int p) \cdot \prod_{i \in I} x_i  \, $ lies in $ M + L_A $,
where $I$ is any index set as above.

\noindent {\bf Algorithm 1. } {\sl (Computing the integral
of a cohomology class of $X_w$)} \hfill \break
\noindent {\sl Input: } A polynomial $p(x_1,\ldots,x_n)$
with  coefficients in a field $k \supset \Q$  \hfill \break
\noindent {\sl Output: }
The integral $\, \int_{X_w} p \,$ of
the corresponding cohomology class on $X_w$.
\begin{enumerate}
\item Compute any Gr\"obner basis ${\cal G}$ for the ideal $M + L_A$.
\item Let $m$ denote the unique standard monomial of degree $n-d$.
\item Find any minimal prime
 $\langle x_j \, : \, j \not\in I \rangle \,$ of $M$,
and compute the normal form of $\,\prod_{i \in I } x_i\,$
modulo the Gr\"obner basis ${\cal G}$. It looks like
$\, \gamma \cdot m$, where $\gamma$ is a non-zero element of $k$.
\item Compute the normal form of $p$ modulo the
Gr\"obner basis ${\cal G}$, and let $\delta \in k $ be the
coefficient of $m$ in that normal form.
\item Output the scalar $\, \delta / \gamma \,\in \, k$.
\end{enumerate}
\vskip .2cm

\noindent

To compute the number of lattice points in $P_w$,
we note that there is a special element in the cohomology
ring $\, H^*( X_w; k) $, denoted $\,td(x_1,\ldots,x_n) \,$
and called the {\it Todd class} of the toric variety $X_w$.  The Todd
class is represented (non-uniquely) by a (non-homogeneous) polynomial
with rational coefficients in the variables $x_1,\ldots,x_n$.  The
polynomial $\,td(x_1,\ldots,x_n)\,$ does  what we want:
$$ \phi_A(w_1,\ldots,w_n) \,\, = \,\,
\# (P_w \,\cap\, \Z^n)
\quad = \quad
\int_{X_w} td(x_1,\ldots,x_n) \cdot exp \, (\sum_{i=1}^n w_i x_i)
 $$
Here the {\it exponential} of a linear form in
(\ref{(4)}) is defined by
the terminating series
\begin{equation}
\label{(7)}
 exp \,(\sum_{i=1}^n w_i x_i) \quad = \quad
\sum_{r=0}^{n-d} {1 \over r \, ! } \cdot
( w_1 x_1 + w_2 x_2 + \cdots +  w_n x_n)^r.
\end{equation}
Pommersheim \cite{BarviPom} gives an  algorithm for computing the
Todd class, which works efficiently even for non-unimodular $A$.
For our applications, however, we prefer
to use the basic formula given in the first line
on page 110 in Fulton's book \cite{Ful}:
\begin{equation}
\label{(8)}
 td(x_1,\ldots,x_n) \,\,\, = \,\,\,
\prod_{i=1}^n { x_i \over 1 - exp(-x_i) }
\,\,\, = \,\,\,
\prod_{i=1}^n \bigl(1+ {1 \over 2} x_i
+ {1 \over 12} x_i^2
- { 1 \over 720} x_i^4+
 \cdots \bigr)
\end{equation}
In this expansion we list only terms of degree $\leq n-d$,
so that (\ref{(8)}) becomes a polynomial  in $x_1,\ldots,x_n$
with $\Q$-coefficients. We conclude with our main result.

\begin{thm} \label{thealgorithm}
The following algorithm computes the
polynomial that represents
 $\psi_A$ on a chamber containing
a given non-negative vector $w \in \R^n$: \rm
\begin{enumerate}
\item Determine the linear
inequalities defining the given chamber.
\item Let $M$ be the ideal generated by the
leading monomials of the Gr\"obner basis for $J_A$ with respect to $w$ and
compute the ideal representing the kernel $L_A$ of $A$. Use these
two ideals to construct the cohomology ring.
\item  Apply Algorithm 1 to the product of the polynomials
in (\ref{(7)}) and (\ref{(8)}).
\end{enumerate}
\end{thm}

A main advantage of this algorithm over other methods is
that we can do the computation parametrically, over the field
$k = \Q (w_1,\ldots,w_n)$. Our output is the actual polynomial
for $\psi_A$, not just some numerical evaluation of it.

For our running example we take the polynomial ring
$S = k[A,B,C,D,E]$ over the field $k = {\Q}(a,b,c,d,e)$.
We fix the reverse lexicographic Gr\"obner basis for the ideal $M + L_A$,
where  $L_A = \langle A+B-D, A+C-E \rangle $ and
 $M = in_w(J_A)$ is the monomial ideal in (\ref{(3)}).
The Todd class (\ref{(8)}) is computed from the formula
$$
(1 + A/2 + A^2/12)
(1 + B/2 + B^2/12)
(1 + C/2 + C^2/12)
(1 + D/2 + D^2/12)
(1 + E/2 + E^2/12)
$$
The normal form of this expression with respect to our Gr\"obner basis
equals
\begin{equation}
\label{(8')}
  td(A,B,C,D,E) \quad = \quad
DE + C/2 + D + E/2 + 1.
\end{equation}
Likewise, the exponential of the general divisor
(\ref{(7)}) on our toric surface,
$$ 1 +  (aA+bB+cC+dD+eE) +  {1 \over 2}(aA+bB+cC+dD+eE)^2, $$
has the following normal form with respect to our Gr\"obner basis:
$$
1 \, + \, (a-b+e)E + (b+c-a)C + (b+d) D\, + \,
\bigl(ab+be+ac+cd+de-(a^2-b^2-c^2)/2 \bigr) DE $$
Multiply this expression with (\ref{(8')}), reduce it to normal form,
and extract the coefficient of the standard monomial $DE$.
The result is the desired polynomial that represents
$\psi_A(a,b,c,d,e)$ on the fifth chamber. Now set $d = e = 0$.

\section{Method Two: BBKLP generating functions}

In the BBKLP method we associate with any rational
polyhedron $P$ in $\R^n $ the following rational generating
function in $n$ variables:
$$f(P, x) \quad = \quad \sum_{u \in P \cap \Z^n} x^u $$
where $x^u$ denotes $x_1^{u_1}x_2^{u_2}...x_n^{u_n}$.
Brion \cite{brion88} proved the following result:

\begin{thm} \label{brion}
For any rational polyhedron $P$ in $\R^n$,
$$ f(P,x) \quad = \quad \sum_{v \in {\rm vertices}(P)} f({\rm cone}(P,v),x),$$
where  $\,{\rm cone}(P,v) \,  =  \,\{u \in \R^n : v+\delta u \in P
\ {\rm for \ all \ sufficiently \ small} \ \delta>0\}$.
\end{thm}

Each of the series $\,f({\rm cone}(P,v),x)\,$ is a rational generating
function, which can be computed using commutative algebra methods
(Hilbert series). For us the only relevant case is that of
{\em unimodular} (also called primitive) cones. A unimodular cone
is a pointed simplicial cone with generators $\, \{u_1,\dots,u_k\} \,$
that form a basis for the lattice
$\, \R \{u_1,\dots,u_k\} \cap \Z^n $.
For unimodular cones, the rational generating function takes
the following simple form:
\begin{equation}
\label{conefct}
 f({\rm cone}(P,v),x) \quad = \quad \prod_{i=1}^k \frac{x^v}{1-x^{u_i}}.
\end{equation}
If $P$ is a rational convex polytope then $f(P,x)$ is a polynomial,
and this polynomial has a ``short'' representation as a rational function
by Theorem \ref{brion}.  Evaluating $f(P,x)$ at
$x=(1,1,1,\dots,1)$ gives the number of integer points in
$P$. However, if we are given $f(P,x)$ as a sum of rational
functions as in Theorem \ref{brion}, then this evaluation
is a nontrivial problem since the point $x=(1,1,1,\dots,1)$
is a pole of  (\ref{conefct}).
We present our solution to this problem in
Algorithm 3.

For a one-dimensional example, let $P$ be the line segment $[0,b]$. Then
$$\, f(P,x) \quad = \quad \frac{1}{1-x} \, + \,  \frac{x^b}{1 - x^{-1}}
\quad = \quad 1  + x + x^2 + \cdots + x^b. $$
The value of this polynomial at $x=1$ equals $b+1$, the number
of lattice points in the segment, but the substitution $x \rightarrow 1$
must be performed with care.

Consider now the polytope $\,P_b \, = \, \{ \, x \in \R^n \, : \, x
\geq 0 , \, A x = b \,\} $, where $A$ is unimodular and $b$ is in the
lattice spanned by the columns of $A$ and lies in the relative
interior of a maximal chamber.  Under these hypotheses, $P_b$ is a
simple polytope such that ${\rm cone}(P_b,v)$ is unimodular for
every vertex $v$ of $P_b$. We shall give a combinatorial
formula for the rational functions representing these cones.

Consider any subset $\sigma \subset \{1,\ldots,n\}$
which is a column basis of the matrix $A$, and let
$v_\sigma $ denote the unique vector in $\R^n$ with
support $\sigma$ and satisfying $\, A \cdot v_\sigma = b$.
The entries of $v_\sigma = v_\sigma(b)$ are linear combinations of the
coordinates of $b$ with integer coefficients.
The vertices of the simple polytope $P_b$ are precisely
those vectors $v_\sigma$ which have all coordinates non-negative.
The edges of $P_b$ emanating from a vertex $v_\sigma$
are parallel  to certain  non-zero vectors
with minimal support in the kernel of $A$. These vectors
are called {\it circuits} in  matroid theory.
For any index $\, i \in \{1,\ldots,n \} \backslash \sigma\,$
let $\,C(\sigma, i) \in \{-1,0,+1\}^n \,$
the  associated {\it basic circuit}. This is the
unique  vector in the kernel of $A$ whose support is
a subset of $\,\sigma \cup \{i\}\,$ and whose
$i$-th coordinate is $+1$. The following lemma is straightforward:

\begin{lemma}
For the vertex of $P_b$ indexed by $\sigma$, the
series (\ref{conefct}) equals
\begin{equation}
 \label{LemmaVert}
 f({\rm cone}(P_b,v_\sigma),x) \quad = \quad
x^{v_\sigma (b)} \cdot \prod_{i \not\in \sigma}
\bigl( 1 - x^{C(\sigma, i)} \bigr)^{-1}
\end{equation}
\end{lemma}

In the formula above, only the monomial
$\,x^{v_\sigma (b)} \,$ depends on the
specific right hand side vector $b$.
The other factors depend only on the
chamber which contains $b$. We use the following procedure
for computing the generating function for the set of non-negative
integer solutions to a unimodular system  $\,A  x = b $.

\smallskip

\noindent {\bf Algorithm 2. }  {\rm (Computing the BBKLP
generating function)}

\noindent {\sl Input: }
Unimodular matrix $A$, a representative vector $b$ for a  chamber
of $A$. \hfill \break
\noindent {\sl Output: }
The generating function $f(P_b,x)$ for the set of lattice points in $P_b$.
\begin{enumerate}
\item Compute all the circuits of matrix $A$. This step is
entirely independent of $b$ and can be done a priori, before processing any
particular chambers.
\item List all subsets $\sigma \subset \{1,\ldots,n\}$ which index
 vertices $v_\sigma$ of the polytope $P_b$.
\item For each $\sigma$ in the previous step, compute
the right hand side of (\ref{LemmaVert}).
\item Output the sum
$\,\sum_\sigma f({\rm cone}(P_b,v_\sigma),x) \,$ of these rational functions.
\end{enumerate}

\medskip

We illustrate the output of
this algorithm for the unimodular $3 \times 6$-matrix
$$ A \quad = \quad
\left [\begin {array}{cccccc} 1&1&1&0&0&0\\\noalign{\medskip}-1&0&0&1&
1&0\\\noalign{\medskip}0&-1&0&-1&0&1\end {array}\right ]
$$
and the right hand side vector $b = (b_1,b_2,b_3) $
in the same chamber with $ (1,3,-2)$.
The polytope $P_b$ is three-dimensional and has six vertices.
Their index sets $\sigma$ and corresponding generating functions
$\, f({\rm cone}(P_b,v_\sigma),x) \,$ are listed in Table \ref{tabla}.
The lattice point enumerator $\,f(P_b,x) \,$ is the sum of
these six rational functions. The number of lattice points
in $P_b$ is found to be:
\begin{equation}
\label{answerpoly}
 f \bigl(P_b,(1,1,1,1,1,1) \bigr) \quad = \quad
\frac{1}{6} (b_1 + 2) (b_1 + 1) (2 b_1 + 3 b_2 + 3 b_3  + 3).
\end{equation}

 \begin{table}[h]
 \begin{center}
 \font\ninerm=cmr11
 {\ninerm
 \begin{tabular}{|c|c|}
 \hline
 index set  & rational function \cr
\hline
    {\{2, 4, 5\}}      & ${x_{{2}}}^{b_1}{x_{{4}}}^{-{b_1}-{b_3}}
{x_{{5}}}^{{b_1}+
{b_3}+{b_2}}\left (1-{\frac {x_{{1}}x_{{4}}}{x_{{2}}}}\right )^{
-1}\left (1-{\frac {x_{{3}}x_{{4}}}{x_{{2}}x_{{5}}}}\right )^{-1}
\left (1-{\frac {x_{{4}}x_{{6}}}{x_{{5}}}}\right )^{-1}$  \cr
   {\{1, 4, 6\}}      &  ${x_{{1}}}^{{b_1}}{x_{{4}}}^{{b_1}+{b_2}}
{x_{{6}}}^{{b_1}+{b_3}+{b_2}}\left (1-{\frac {x_{{2}}}{x_{{1}}x_{{4}}}}\right )^{-
1}\left (1-{\frac {x_{{3}}}{x_{{1}}x_{{4}}x_{{6}}}}\right )^{-1}\left
(1-{\frac {x_{{5}}}{x_{{4}}x_{{6}}}}\right )^{-1}$  \cr
    {\{2, 4, 6\}}     & ${x_{{2}}}^{{b_1}}{x_{{4}}}^{{b_2}}
{x_{{6}}}^{{b_1}+{b_3}+{b_2}}\left (1-{\frac {x_{{5}}}{x_{{4}}x_{{6}}}}\right )^{-1}\left (
1-{\frac {x_{{1}}x_{{4}}}{x_{{2}}}}\right )^{-1}\left (1-{\frac {x_{{3
}}}{x_{{2}}x_{{6}}}}\right )^{-1}$  \cr
     {\{3, 4, 5\}}      & $ {x_{{3}}}^{{b_1}}
{x_{{4}}}^{-{b_3}}{x_{{5}}}^{{b_3}+{b_2}}
\left (1-{\frac {x_{{4}}x_{{6}}}{x_{{5}}}}\right )^{-1}\left (1-{
\frac {x_{{2}}x_{{5}}}{x_{{3}}x_{{4}}}}\right )^{-1}\left (1-{\frac {z
_{{1}}x_{{5}}}{x_{{3}}}}\right )^{-1}$  \cr
    {\{3, 4, 6\}}     & ${x_{{3}}}^{{\it b1}}{x_{{4}}}^{{b_2}}
{x_{{6}}}^{{b_3}+{b_2}}
\left (1-{\frac {x_{{2}}x_{{6}}}{x_{{3}}}}\right )^{-1}\left (1-{
\frac {x_{{5}}}{x_{{4}}x_{{6}}}}\right )^{-1}\left (1-{\frac {x_{{1}}z
_{{4}}x_{{6}}}{x_{{3}}}}\right )^{-1}$ \cr
   {\{1, 4, 5\}} & ${x_{{1}}}^{{b_1}}{x_{{4}}}^{-{b_3}}{x_{{5}}}^{{b_1}
+{b_3}+{b_2}}\left (1-{\frac {x_{{3}}}{x_{{1}}x_{{5}}}}\right )^{-1}\left
(1-{\frac {x_{{4}}x_{{6}}}{x_{{5}}}}\right )^{-1}\left (1-{\frac {x_{{
2}}}{x_{{1}}x_{{4}}}}\right )^{-1}$ \cr
  \hline
 \end{tabular}}
 \end{center}
 \caption{Rational functions associated to the
 six supporting cones at the vertices} \label{tabla}
 \end{table}

This last evaluation can be done symbolically, for instance, using the
command {\tt simplify} in {\tt Maple}, but the symbolic simplification
is too slow for larger examples.  We compute the limit $\,x_i
\rightarrow 1 \,$ in such rational functions by first specializing to
a single variable $t$, in a manner to be described in Algorithm 3.

We next discuss how to implement Step 2 of the algorithm, namely,
how to efficiently list all vertices of $P_b$. The first possibility
is to compute the prime decomposition of the monomial ideal $M$
which was used in Section 2 to encode the chamber of $b$. Indeed,
a subset $\sigma$ corresponds to a vertex of $P_b$ if and only if
$\,\langle x_i \, : \, i \not\in \sigma \rangle \,$ is a minimal
prime of $M$. The second possibility is to precompute the
vector-valued linear functions $\,v_\sigma(b) \,$ for all
column bases of $A$. Similarly to the computation of the circuits
in step 1, this can be done a priori, before processing any
particular chambers. For any particular chamber, we
take the sum in step 4 only over those bases $\sigma$ which
satisfy $v_\sigma(b) \geq 0 $. In our practical implementation
we opted for a third possibility, namely, to apply
a depth-first search algorithm to the graph
of  basic feasible solutions of $\,Ax = b, x \geq 0 $, where
the edges are basis exchanges as in the simplex algorithm
 \cite[Chapter 8]{schrijver}. A considerable speed-up
over our crude {\tt Maple} implementation
can still be obtained by using the reverse-search algorithm
of Avis and Fukuda \cite{avisfukuda}.

The output of Algorithm 2 is a generating function $\,f(P_b,x) \,$
that represents the vector partition function on a particular
chamber. Now we face the  problem  to evaluate,
for any particular $b \in \Z^d$, the
limit of $f(P_b,x)$ as $x$ tends to $(1,1,\dots,1)$.
In the literature there are two approaches to this problem:
the  Barvinok-Brion method \cite[Algorithm 5.2]{BarviPom}
and the Dyer-Kannan method \cite{dyerkannan}.
Both methods consider the rational series as a sum of
exponential functions each of which converges for almost all
choices of $x$. The first
approach essentially takes the residue of the function and the second
computes the value of the rational function at a point close to
$x=(1,1,\dots,1)$ and carefully rounds the answer to the nearest
integer.  When we tried these two approaches experimentally,
we ran into memory problems and numerical instabilities.
In our experience, the following alternative method
works rather well in practice:

\noindent {\bf Algorithm 3.}  {\rm (Evaluating
the BBKLP generating function at $(1,1,1,\ldots,1)$)}
\begin{enumerate}
\item Eliminate  $\,rank(A)$ many variables
by substitutions $x_i=1$ where $i$ runs over a column
basis of $A$. All denominators
$\, 1 - x^{C(\sigma, i)} \,$ remain nonzero.

\item For each vertex $v_\sigma$ of $P_b$, replace
each remaining variable $x_j$ by $\, 1 - j  \,t $.
This transforms (\ref{LemmaVert}) into
a rational function in one variable $t$.
We express the result in the
form numerator/denominator, where the numerator and denominator are
relatively prime polynomials in $t$ with integer coefficients.

\item Replace the sum of rational functions, one
for each vertex of $P_b$,
by a single rational function $p(t)/q(t)$. Here $q(t)$ is the
least common multiple of the denominators of the
rational functions produced in step 2.

\item  Both $p(t)$ and $q(t)$ vanish at $t=0$.
Let $t^\alpha$ be the largest common factor. We compute the limit of
$p(t)/q(t)$ as $t \rightarrow 0$ using L'H\^opital's rule.
For that we need the value at $0$ of the $\alpha$-th
derivatives of $p$ and $q$. These can be found in {\tt Maple}
using the built-in feature of {\it automatic differentiation}.
This allows us to retain the representation of $p$ as   sum of
terms, one for each vertex of $P_b$, and that of
$q$ as a product of binomials $\, 1 - x^{C(\sigma, i)} $.

\item Output $(f/g)(0)$.
\end{enumerate}

At the beginning of this section we had assumed that $b$
lies in the relative interior of a maximal chamber.
This assumption  can be removed easily.
It was  made in order to uniquely identify the chamber and
hence a  representation of $P_b$ as a simple polytope.
 If $b$ happens to lie in
a lower-dimensional chamber, and $P_b$ is
not simple, then we can use the combinatorial description
of any adjacent maximal chamber in step 2 of Algorithm 2.
This is consistent with the fact,
implied by Theorem \ref{structure},
that the polynomials
representing $\phi_A(b)$ on different chambers must
agree on the intersection of the closures of these chambers.

\section{Contingency Tables}

In the remainder of this paper, we report on the implementation
and performance of our methods  for two important
families of unimodular matrices $A$. We present
both computational and mathematical results.
We ran all our experiments in a computer with a single
Pentium-III CPU with 700Mhz
and 256 MB RAM using the computer algebra packages {\tt Macaulay 2}
and {\tt Maple}. The generation of chambers was performed
using {\tt Topcom} and {\tt Puntos}; see \cite{puntos, rambau}.

Let $r=(r\sb 1,\ldots, r\sb m)$ and $c=(c\sb 1,\ldots, c\sb n)$ be
compositions of a fixed integer $N \geq 1$. Let $\Sigma\sb
{rc}$ denote the set of all $m\times n$ non-negative integer matrices
in which row $i$ has sum $r\sb i$ and column $j$ has sum $c\sb
j$. Thus $\sum\sb {i,j}T\sb {ij}=N$ for any $T\in\Sigma\sb {rc}$.
We are interested in the number $\# \Sigma\sb {rc} $ of
matrices in $\Sigma\sb {rc}$. This number equals $\phi_A(b)$  where
$A$ is the node-edge incidence matrix of the complete bipartite graph
$K_{n,m}$ and the vector $b$ is the vector $(r,c)$. Thus we are
counting the lattice points in a {\it transportation polytope}.
There is an extensive literature on
computing the function $(r,c) \mapsto \# \Sigma\sb {r,c} $. See
\cite{persi2} and the  references therein.

We implemented the Gr\"obner bases algorithm described in Section 2
in the computer algebra system {\tt Macaulay 2},
which was developed by Grayson and Stillman
\cite{macaulay}.  Our {\tt Macaulay 2} program for computing
the polynomial representing $\phi_A$ on a single chamber is very short
and simple. In Appendix 2 we list the entire program for one chamber in
the $4 \times 4$-contingency table case.

As mentioned in the introduction, $4 \times 4$-tables are an important
benchmark. There are $3694$ chambers modulo symmetry. On each chamber,
the function $\,(r,c) \mapsto \# \Sigma_{rc}\,$ is a polynomial
of degree nine in the eight variables $r_1,r_2,r_3,r_4$, $c_1,c_2,c_3,c_4$.
Mount \cite{mount2} computed (interpolation schemes for)
all $3694$ polynomials. He reported
a 3 hour calculation for each chamber, adding up to a total of
6 weeks of  distributed computing for preprocessing all chamber
polynomials.

Our experiments show that the
Gr\"obner basis computation is as least
as fast as  Mount's interpolation technique.
We computed all $3694$ chamber polynomials using
the {\tt Macaulay 2} code listed in Appendix 2.
The running time per chamber ranged from
7 seconds to 45 minutes.  It took us  6 1/2 weeks
sequential computing time to complete the task.
Our {\tt Macaulay 2} code can easily be modified to get
the numerical value $\# \Sigma\sb {rc} $ for any given
$r,c \in \N^4$.  Computing such
numerical instances takes 20 seconds on the average
for $4 \times 4$-tables. Similar computations for
$4 \times 5$-tables have not yet been successful
in {\tt Macaulay 2}.

\smallskip

We implemented the BBKLP method described in Section 3 for contingency
tables in {\tt Maple}. The generation of the BBKLP rational function
(Algorithm 2) runs rather well for our purpose. It takes only a
few seconds for $4 \times 4$-tables, as little as five minutes for $4
\times 5$-tables and up to two days for $5 \times 5$-tables.
We wish to stress that our {\tt Maple} code does not use
optimal techniques for vertex enumeration of polytopes.
For instance, using the Avis-Fukuda reverse search algorithm
 \cite{avisfukuda}
instead of depth-first search would give a significant speed-up
over our crude implementation. For example, the
vertices of a $5 \times 5$ transportation polytope can be computed in
a few seconds using the program {\tt lrs} \cite{avis}.

The second stage in the BBKLP method is Algorithm 3.  This can be
applied either for {\sl symbolic parameters} $r_i$ and $c_j$, in which
case the output is a chamber polynomial, or for {\sl numerical values}
of $c_i$ and $r_j$, in which case the output is the integer $\,\#
\Sigma_{r,c} $.  The second application of Algorithm 3 performs
extremely well in {\tt Maple}. The running time of a numerical
evaluation $\,(r,c) \mapsto \# \Sigma_{r,c}\,$ using Algorithm 3 is
close to one minute for $4 \times 4$ tables, about ten minutes for $4
\times 5$ tables, and about ten days for $5 \times 5$-tables. On the
other hand, the first (symbolic) application of Algorithm 3 is only
possible for smaller matrices, and is generally outperformed by the
Gr\"obner basis computation in {\tt Macaulay 2}.

Here are three test cases that show the power of the BBKLP technique,
with numerical evaluation in Algorithm 3. The largest instance
computed by Mount \cite{mount2} is the number of $4\times 5$-tables
with margins $[3046, 5173,6116,10928]$ and $[182,778,3635,9558,11110]$.
It took him 20 minutes of parallel computing to find the value
$\, 23196436596128897574829611531938753$.
Our {\tt Maple} program reproduces this number in only 10 minutes.

Consider next the  $4 \times 5$-tables whose margins are
$[338106,$ $ 574203, 678876,$ $ 1213008]$
and $[20202,142746, 410755, 1007773, 1222717]$.
Their number equals
$$ 316052820930116909459822049052149787748004963058022997262397.$$
The computation took 35
minutes. The associated transportation polytope $P_b$ is
$12$-dimensional and has $ 976$ vertices.

Finally, we counted all $5 \times 5$-tables with margins
$[30201,59791,70017,41731,$ $58270]$ and
$[81016,68993,47000,43001,20000]$.  The associated 16-dimensional
transportation polytope has 13150 vertices. This computation took 10
days and the answer is a 64 digit number. Algorithm 2 ran about 2 1/2
days.  The size of its output exceeds the memory of our computer.
Therefore we had to apply the lcm-computation in Algorithm 3 to
incremental pieces of this output.

Our {\tt Maple} program for counting $4 \times 4$ and $4 \times 5$-tables
is available at \url{www.math.ucdavis.edu/~deloera/contingency.html}. This webpage includes all relevant data
for the two specific $4 \times 5$-tables discussed above.

\smallskip

The subproblem of enumerating all chambers lead us to take
a look at the structure of the chamber complex for the $n \times m$
contingency tables. This chamber complex is the cone over the
chamber complex of the product of two simplices 
$$ \Delta_{n-1} \times \Delta_{m-1} \quad = \quad
\bigl\{ \, (x_1,\ldots, x_n ; y_1,\ldots,y_m) \in \R_{\geq 0}^{n+m} 
\,\,: \,\, \sum_{i=1}^n x_i =  \sum_{j=1}^m y_j = 1 \, \bigr \} . $$
The combinatorial structure of the polytope $\Delta_{n-1} \times \Delta_{m-1}$
can be read  off from the complete bipartite graph $K_{n,m}$. For instance, 
the full-dimensional  simplices
in $\Delta_{n-1} \times \Delta_{m-1}$ correspond
to spanning trees of $K_{n,m}$, while the facets of $\Delta_{n-1} \times
\Delta_{m-1}$ are complete bipartite subgraphs of $K_{n,m}$ obtained
by removing a vertex.  The $(n+m-3)$-dimensional subsimplices
correspond to a spanning tree minus an edge. 
We define a {\em diagonal section} of $\Delta_{n-1} \times \Delta_{m-1}$
to be any affine hyperplane  which is spanned by vertices of
$\Delta_{n-1} \times \Delta_{m-1}$ but is not a facet hyperplane.
The diagonal sections are in bijection with 
spanning forests of $K_{n,m}$ which have exactly two components.
Let $\Omega_{n,m}$ denote the subdivision of the polytope
$\Delta_{n-1} \times \Delta_{m-1}$ defined by the diagonal sections.
Equivalently, two points $(x,y)$ and $(x',y')$ in
$\Delta_{n-1} \times \Delta_{m-1}$ lie in the same open cell
of $\Omega_{n,m}$ if and only if the lie on the same
side of any hyperplane of the form
$$  x_{i_1} + \cdots + x_{i_r} \, + \,
 y_{j_1} + \cdots + y_{j_s} \, = \, 0 .$$
We call $\,\Omega_{n,m}\,$ the {\it diagonal section
complex} of $\, \Delta_{n} \times \Delta_{m}$.

\begin{prop} The chamber complex of 
 $\Delta_{n} \times \Delta_{m}$ coincides with the
diagonal section complex $\Omega_{n,m}$.
There exist virtual chambers  whenever $\,m+n \geq 7 $.
\end{prop}

\proof For any polytope whatsoever, the
diagonal section complex can be defined,
and it always refines the chamber complex.
The two complexes are equal for polygons,
but they are usually not equal for
higher dimensional polytopes.
What we are claiming is that they are equal 
for products of two simplices.

The key observation is this: the intersection of a
diagonal section with any facet of the polytope
$\Delta_{n} \times \Delta_{m}$ equals the
convex hull of all vertices of the facet
which lie in that diagonal section.
This follows from our graph-theoretical dictionary,
since each facet corresponds to a complete
bipartite subgraph $K_{n-1,m}$ or $K_{n,m-1}$.
From this it follows that each codimension one
simplex spanned by vertices of  $\Delta_{n} \times \Delta_{m}$ 
has the same intersection with 
the boundary of $\Delta_{n} \times \Delta_{m}$ 
as the corresponding diagonal section does.
Therefore the chamber complex equals the diagonal section complex.
The assertion about virtual chambers is proved by
computer calculations for $K_{2,5}$ and $K_{3,4}$. \qed

\section{Kostant's Partition Function}

Let $A$ be the node-arc incidence matrix of the complete acyclic graph
$K_n$.  The function $ \phi_{K_n} := \phi_A $ is the {\it Kostant
partition function} for the root system $A_{n-1}$. Explicitly, let
$e_1,\ldots,e_n$ denote the standard basis of $\Z^n$, and let $
E_{1,2}, E_{1,3}, \ldots, E_{n,n-1}$ denote the standard basis of
$\Z^{{n \choose 2}}$. The matrix $A$ represents the map
$$ \tau \quad : \quad \R_{\geq 0}^{{n \choose 2}} \rightarrow \R^n,
\, E_{i,j} \mapsto e_i - e_j \, \,\, \hbox{for} \,\, 1 \leq i < j \leq
n . $$
The image of $\tau$ is the $(n-1)$-dimensional cone
$$ S_n \,\, = \,\,
\bigl\{ (u_1,\ldots,u_n) \in \R^n \,|\,
u_1 + \cdots + u_i \geq 0
\,\,\,\hbox{for}\,\, 1 \leq i < n \,\,\,
\hbox{and}\,\,
u_1 + \cdots + u_n = 0 \bigr\}. $$

Kirillov \cite[page 57]{kirillov} posed the
problem of finding the number of chambers for Kostant's partition function.
 We give a partial solution to Kirillov's problem by
determining the number of chambers for $n \leq 7$.
See Table \ref{kiritable} below.

We also
computed all chamber polynomials representing $\phi_{K_n}(b)$ for
$ n \leq 6$. This was done using our {\tt Macaulay 2} implementation
(see Appendix 2) of the Gr\"obner basis method in Section 2. For instance, for
$n=6$, there are $820$ chamber polynomials, each of degree $10$
in five variables.
All of these polynomials are available, both in expanded form
and as an on-line calculator, at our web site
\url{www.math.ucdavis.edu/~deloera/kostant.html}.

As a small sample of our results we present all
chambers and chamber polynomials for $n=4$.
These polynomials were first computed by mathematical physicists
in \cite{schmidtbincer}. Analogous computations
for $n \geq 5$ had been infeasible in  1984. In Appendix 1
we list all those chamber polynomials for $n=5$ which can be
factored over $\Q$. Several authors \cite{kirillov,schmidtbincer}
have studied factorization patterns of polynomials representing
Kostant's partition function. A forthcoming paper by Postnikov and Stanley
contains the state of the art. Our data provide complementary information
to their  combinatorial results.

\begin{figure}
 \begin{center}
     \includegraphics[scale=.5]{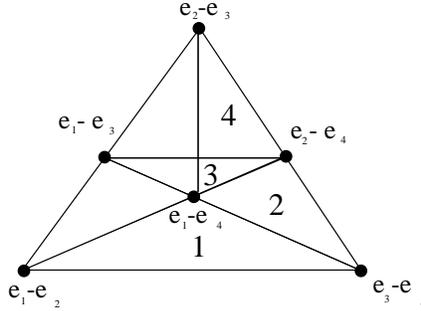}
 \caption{The chamber complex for the 
complete graph $K_4$} \label{A4chambers}
 \end{center}
 \end{figure}

The cone $S_4$ spanned by the columns of the node-arc incidence matrix
of $K_4$ is a three-dimensional triangular cone. The chamber complex
is a subdivision of this cone into seven triangular cones. See Figure
\ref{A4chambers} for a 2-dimensional perpendicular slice showing the
chamber complex.  The formulas below are given only in terms of $b_1,
b_2, b_3$, in view of $b_4=- b_1-b_2-b_3$. By the symmetry of the
example it is enough to give the four polynomials for the indicated
chambers in Figure \ref{A4chambers}. The label of a chamber in the
figure and its polynomial match.

 {\small
 \begin{enumerate}
 \item If $\,{\rm min} \{b_3, -b_2, b_1 + b_2 \}
 \geq 0 \,$ then $ \phi_{K_4}(b) \,\, = \,\,
 (b_1+b_2+3)(b_1+b_2+2)(b_1+b_2+1)/6. $

 \item If $\,{\rm min} \{b_1, b_2,b_3 \} \geq 0 \,$ then
 $ \phi_{K_4}(b) \quad = \quad
 (b_1+1)(b_1+2)(b_1+3b_2+3)/6 $

 \item If $\,{\rm min} \{b_1, b_2, b_1 + b_3, b_2 + b_3,- b_3 \} \geq
 0 \,$ then $ \phi_{K_4}(b) \quad = \quad
1+{\frac {11}{6}}\,{b_1}+2/3\,{b_3}+{b_2}+3/2\,{b_1}\,{
b_2}+{{b_1}}^{2}+1/6\,{{b_1}}^{3}+1/2\,{{b_1}}^{2}{b_2}
-1/6\,{{b_3}}^{3}$
$-1/2\,{b_1}\,{{b_3}}^{2}+1/2\,{b_1}\,{
b_3}-1/2\,{{b_3}}^{2}$

 \item If $\,{\rm min} \{
 b_1, b_2 + b_3, - b_1 - b_3 \}
 \geq 0 \,$ then $ \phi_{K_4}(b) \, = \,
 (b_1 + 2) (b_1 + 1) (2 b_1 + 3 b_2 + 3 + 3 b_3)$
 \end{enumerate}}

Let $\Gamma(K_n)$ be the chamber complex for $\phi_{K_n}$.
This is a polyhedral decomposition of
the cone $S_n$. We have the following result:

 \begin{thm}
\label{Kirikiri}
The complex
 $\Gamma(K_n)$ has chambers with at least $2^{\lfloor n/2 \rfloor}$
 facets. There exist virtual chambers  for $n \geq 5$.
 The exact number of chambers
   for $n \leq 7$ is given by Table \ref{kiritable}.

 \begin{table}[h] 
 \begin{center}
 \font\ninerm=cmr11
 {\ninerm
 \begin{tabular}{|c|c|c|} 
 \hline
  {\bf $n$}   & Number of chambers & Degree of $\phi_{K_n}$ \cr
 \hline
    3      & 2       & 1 \cr
    4      & 7       & 3 \cr
    5      & 48      & 6\cr
    6      & 820     & 10 \cr
    7      & 44288   & 15 \cr
  \hline
 \end{tabular}}
 \caption{Chambers for $A_n$} \label{kiritable}
 \end{center}
 \end{table}
 \end{thm}

\begin{figure}
\begin{center}
    \includegraphics[scale=.3]{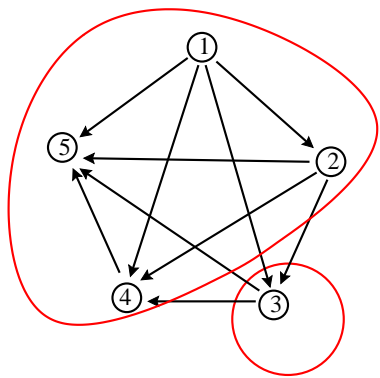}
\caption{A balanced cut for $K_5$} \label{balanced}
\end{center}
\end{figure}

\proof Let $A_{n-1} = \{\, e_i - e_j \,:\, 1 \leq i < j \leq n \, \} $.
There is a
well-known bijection between cuts of the digraph $K_n$ and hyperplanes
spanned by subsets of $A_{n-1}$. For odd values of $n$ there
are ``balanced'' cuts for $K_n$. By a balanced cut we mean one where
the hyperplane associated divides the set of roots $e_i - e_j$ outside
the hyperplane into equal size groups. In Figure \ref{balanced} we
show one such cut for $K_5$ that leaves two roots in each side of the
plane $x_3=0$. To obtain such a balanced cut for general $K_n$, odd
$n$, note that there is a middle node labeled $\lfloor n/2 \rfloor +1$
that has exactly as many entering arcs as leaving arcs. The cut
$\{1,2,3,\dots,\lfloor n/2 \rfloor,\lfloor n/2 \rfloor+2,\dots,n\}$
and $\{ \lfloor n/2 \rfloor +1\}$ is balanced. The vectors in
$A_{n-1}$ that lie on the plane $x_{\lfloor n/2 \rfloor +1}=0$ 
form the configuration $A_{n-2}$.

The intersection of the chamber complex of $A_{n-1}$ with a balanced
hyperplane $H$ induces exactly the chamber complex of $A_{n-2}$.
Indeed, the only way to create new cells for $H \cap A_{n-1}$ (not
already in $A_{n-2}$) is if simplices with vertices on opposite
halfspaces of $H$ cut out new vertices in $H$. But pairs of vectors on
opposite sides of $H$ are always collinear with a root $e_i - e_j$
lying on $H$.  The collinearities can be read off from cycles of
length three in the graph $K_n$ that touch the vertex $\lfloor n/2
\rfloor +1$. The existence of triples of collinear vectors, the
center one inside $H$, has another effect: a chamber $\gamma$ of $H$,
one of whose vertices is part of a collinearity, extends to both
halfspaces of $H$. This is because the $(n-2)$-simplices inside $H$
that make up that chamber can be turned into $(n-1)$-dimensional
simplices by coning them with the two extremes of the collinearity
that do not belong to the hyperplane $H$.  Note that the completion
happens in both halfspaces of $H$ but the result of intersecting these
simplices has in common the open cell $\gamma$ that connects both
sides. This might not be the final chamber that extends $\gamma$, as
other vectors in $A_{n-1}$ not lying on $H$ could be used to build and
intersect more $(n-1)$ simplices, but the the result will be contained
in this initial convex cell that touches both halfspaces of $H$. The
number of facets will be then at least twice the number of facets of
$\gamma$. The doubling on the number of facets occurs for odd values
of $n$ but for even values at worse remains the same. Thus,
recursively we can build a chamber with exponentially many facets.
The rest of the statement follows from computer calculations based on
the duality between chambers and triangulations as explained in
Section 2.\qed

\bigskip

\noindent{\bf Acknowledgements:} We thank A. Barvinok, R. Hemmecke,
R.P. Stanley, M.~Stillman, D. Zeilberger, and G. Ziegler for helpful
conversations. Jesus De Loera was supported by NSF Grant DMS-0073815.
Bernd Sturmfels was supported by NSF Grant DMS-9970254 and the Miller
Institute at UC Berkeley.

\section{Appendix: Kostant partition function for $A_4$}

Here we consider $\phi_{K_n}$ for $n = 5 $.  This is Kostant's
partition function for the root system $A_4$. The chamber complex can
be visualized as a subdivision of a tetrahedron. This polyhedral
complex has $19$ vertices, $77$ edges, $107$ triangles and $48$
three-dimensional chambers.  Only two of these $48$ chambers are not
tetrahedra: they are bipyramids.  Thirty of the $48$ chamber
polynomials are irreducible over $\Q$. We explicitly list the other
$18$ chamber polynomials, namely those that factor, together with
defining inequalities for their chambers.

{\small
\begin{enumerate}
\item If $ {\rm min} \left \{{b_1},{b_2},{b_3},{b_4}\right \}\geq 0  $ then

${\frac {1}{360}} \left ({b_1}+3\right )\left ({b_1}+2\right )
\left ({b_1}+1\right )\left ({{b_1}}^{2}+5 {b_1} {b_2}+9
 {b_1}+20+10 b_2^{2}+30 {b_2}\right )$

$\left ({b_2}+3
+{b_1}+3 {b_3}\right )$

\item If $ {\rm min} \left \{{b_4},{b_1}+{b_2}+{b_3},-{b_3},-{b_2}\right
\} \geq 0$ then

${\frac {1}{360}} \left ({b_2}+5+{b_1}+{b_3}\right )\left ({
b_2}+4+{b_1}+{b_3}\right )\left ({b_2}+3+{b_1}+{b_3}
\right ) ({b_2}+2+{u1}+ {b_3})$

$\left ({b_2}+1+{
b_1}+{b_3}\right )\left ({b_2}+3+{b_1}-2 {b_3}\right )$


\item If $ {\rm min} \left \{-{b_3},-{b_2}-{b_4},-{b_1}-{b_4},{b_1}+{\it
b_2}+{b_3}+{b_4}\right \}\geq 0  $ then

${\frac {1}{360}} \left ({b_4}+3+{b_1}+{b_3}+{b_2}\right )
\left ({b_4}+2+{b_1}+{b_3}+{b_2}\right )({b_4}+1+
{b_1}+{b_3}+{b_2})$

$ (60+56 {b_1}+6 {b_2}-
14 {b_3}-54 {b_4}+9 {b_3} {b_4} {b_1}+6 {b_2}
 {b_3}-3 {b_2} {b_4} {b_1}-9 b_3^{2}{\it
b_1}+$

$3 {b_1} b_4^{2}+3 {b_2} b_4^{2}-6 {b_3}
b_4^{2}+27 {b_1} {b_2}-9 {b_1} {b_3}-9 b_1^{2}{b_4}
+6 b_2^{2}{b_4}+9 b_1^{2}{b_2}-6 {
b_1} b_2^{2}+3 b_3^{2}{b_4}+6 b_3^{2}{
b_2}+24 {b_3} {b_4}-45 {b_1} {b_4}-6 {b_2} {
b_3} {b_4}-b_3^{3}+6 b_1^{2}-9 b_2^{2}-
15 b_3^{2}-2 b_1^{3}+3 b_2^{3}+9 b_4^{
2} )$

\item If $ {\rm min}\left \{-{b_3},-{b_2},-{b_4},{b_1}+{b_2}+{b_3}+{\it
b_4}\right \} \geq 0  $ then

${\frac {1}{360}} ({b_2}+3+{b_1}-2 {b_3} )\left
({b_4}+3+{b_1}+{b_3}+{b_2}\right )\left ({b_4}+2+{\it
b_1}+{b_3}+{b_2}\right )$

$({b_4}+1+{b_1}+{b_3}+{\it
b_2}) (b_2^{2}+2 {b_2} {b_3}+9 {b_2}+2
{b_1} {b_2}-3 {b_2} {b_4}+9 {b_3}+b_1^{2}+20
+9 {b_1}+2 {b_1} {b_3}+b_3^{2}-3 {b_3} {b_4
}-21 {b_4}-3 {b_1} {b_4}+6 b_4^{2})$


\item If $ {\rm min} \left \{{b_3},{b_4},-{b_2},{b_1}+{b_2}\right \} \geq 0  $ then

${\frac {1}{360}} \left ({b_2}+{b_1}+3\right )\left ({b_2}+2+
{b_1}\right )\left ({b_2}+1+{b_1}\right )\left ({b_2}+5+{
b_1}\right )\left ({b_2}+4+{b_1}\right )$

$\left ({b_2}+3+{
b_1}+3 {b_3}\right )$


\item If $ {\rm min} \left \{{b_1},{b_3},{b_2}+{b_4},-{b_1}-{b_3}-{b_4
}\right \} \geq 0  $ then

${\frac {1}{360}} \left ({b_1}+3\right )\left ({b_1}+2\right )
\left ({b_1}+1\right )$

$ (60+56 {b_1}+110 {b_2}+70 {
b_3}+50 {b_4}-30 {b_3} {b_4} {b_1}+90 {b_2} {
b_3}+30 b_2^{2}{b_3}-15 {b_2} {b_4} {\it
b_1}$

$-6
 b_1^{2}{b_3}-15 b_3^{2}{b_1}-30 {b_1} {{
b_4}}^{2}-30 {b_2} b_4^{2}-30 {b_3} b_4^{2}
+57 {b_1} {b_2}+21 {b_1} {b_3}-15 b_1^{2}{\it
b_4}+3 b_1^{2}{b_2}+15 {b_1} b_2^{2}-30 {{\it
b_3}}^{2}{b_4}+30 {b_2} {b_4}-15 {b_1} {b_4}+15 {
b_2} {b_3} {b_1}-10 b_3^{3}-20 b_4^{3}+6 {
{b_1}}^{2}+60 b_2^{2}-2 b_1^{3}+10 b_2^{3}-
30 b_4^{2} )$


\item If $ {\rm min} \left \{{b_1},{b_3}+{b_4},{b_2}+{b_4},-{b_1}-{b_4
}\right \} \geq 0  $ then

${\frac {1}{360}} \left ({b_1}+3\right )\left ({b_1}+2\right )
\left ({b_1}+1\right )$

$ (3 b_1^{2}{b_3}-6 {{b_1}
}^{2}{b_4}+9 b_1^{2}+3 b_1^{2}{b_2}+51 {b_1}
+57 {b_1} {b_2}+15 {b_1} b_2^{2}-9 {b_1} {
b_4}$

$-15 {b_1} b_4^{2}+27 {b_1} {b_3}-15 {\it
b_2} {b_4} {b_1}+15 {b_2} {b_3} {b_1}+30 {b_2}
 {b_4}+60+60 {b_3}+60 b_2^{2}+40 {b_4}+90 {b_2
} {b_3}-30 b_4^{2}+110 {b_2}-10 b_4^{3}+30 {{
b_2}}^{2}{b_3}+10 b_2^{3}-30 {b_2} b_4^{2}
)$


\item If $ {\rm min} \left \{-{b_2}-{b_4},-{b_2}-{b_3},-{b_1}-{b_3}-{\it
b_4},{b_1}+{b_2}+{b_3}+{b_4}\right \} \geq 0  $ then

${\frac {1}{360}} \left ({b_4}+3+{b_1}+{b_3}+{b_2}\right )
\left ({b_4}+2+{b_1}+{b_3}+{b_2}\right )\left ({b_4}+1+
{b_1}+{b_3}+{b_2}\right )$

$ (60+51 {b_1}+11 {b_2}-
9 {b_3}-59 {b_4}+3 {b_3} {b_4} {b_1}-6 b_2
^{2}{b_3}-9 {b_2} {b_4} {b_1}-3 b_1^{2}{b_3}
$

$-6 b_3^{2}{b_1}+9 {b_1} b_4^{2}-9 {b_3} {{
b_4}}^{2}+27 {b_1} {b_2}-9 {b_1} {b_3}-3 {{b_1}
}^{2}{b_4}+9 b_2^{2}{b_4}+6 b_1^{2}{b_2}-3 {
b_1} b_2^{2}+6 b_3^{2}{b_4}+24 {b_3} {\it
b_4}-39 {b_1} {b_4}+6 {b_2} {b_3} {b_1}-3 {{b_3
}}^{3}+2 b_4^{3}+9 b_1^{2}-12 b_2^{2}-18 {{
b_3}}^{2}+b_2^{3}+12 b_4^{2})$


\item If $ {\rm min} \left \{-{b_2},-{b_3}-{b_4},{b_1}+{b_2}+{b_3}+{\it
b_4},-{b_1}-{b_2}-{b_4}\right \} \geq 0  $ then

${\frac {1}{360}} \left ({b_4}+3+{b_1}+{b_3}+{b_2}\right )
\left ({b_4}+2+{b_1}+{b_3}+{b_2}\right )\left ({b_4}+1+
{b_1}+{b_3}+{b_2}\right )$

$\left ({b_2}+4-{b_4}-{b_3}+ {b_1}\right )
 (2 b_2^{2}+13 {b_2}+4 {b_1} {
b_2}-{b_2} {b_4}-{b_2} {b_3}+15+13 {b_1}$

$-{b_1}
 {b_3}+2 b_1^{2}-7 {b_3}-7 {b_4}-{b_1} {b_4
}+2 b_4^{2}+4 {b_3} {b_4}+2 b_3^{2} )$


\item If $ {\rm min}\left \{{b_4},{b_1}+{b_2}+{b_3},-{b_2}-{b_3},-{\it
b_1}-{b_3}\right \} \geq 0  $ then

$-{\frac {1}{360}} \left ({b_2}-3+{b_3}-2 {b_1}\right )
\left ({b_2}+5+{b_1}+{b_3}\right )\left ({b_2}+4+{b_1}+
{b_3}\right )$

$\left ({b_2}+3+{b_1}+{b_3}\right )\left ({
b_2}+2+{b_1}+{b_3}\right )$
$\left ({b_2}+1+{b_1}+{b_3}
\right )$


\item If $ {\rm min} \left \{{b_1},{b_2}+{b_3}+{b_4},-{b_1}-{b_4},-{\it
b_1}-{b_3}\right \} \geq 0  $ then

${\frac {1}{360}} \left ({b_1}+3\right )\left ({b_1}+2\right )
\left ({b_1}+1\right )\left ({b_3}-2 {b_4}+{b_2}+3\right
)$
$ (10 b_2^{2}+30 {b_2}+15 {b_1} {b_2}+20 {
b_2} {b_3}+20 {b_2} {b_4}+20+15 {b_1} {b_4}+30
 {b_3}+30 {b_4}+15 {b_1} {b_3}+6 b_1^{2}+10 {
{b_4}}^{2}+10 b_3^{2}+24 {b_1}+20 {b_3} {b_4}
)$


\item If $ {\rm min} \left \{{b_1},{b_4},{b_2}+{b_3},-{b_1}-{b_3}\right
\} \geq 0  $ then

${\frac {1}{360}} \left ({b_1}+3\right )\left ({b_1}+2\right )
\left ({b_1}+1\right )\left ({b_1}+4+2 {b_2}+2 {b_3}
\right )$

$\left (5 b_2^{2}+5 {b_1} {b_2}+10 {b_2} {
b_3}+20 {b_2}+13 {b_1}+5 {b_1} {b_3}+20 {b_3}+2
 b_1^{2}+5 b_3^{2}+15\right )$


\item If $ {\rm min} \left \{-{b_2}-{b_3},-{b_4},{b_1}+{b_2}+{b_3}+{\it
b_4},-{b_1}-{b_3}\right \} \geq 0  $ then

$-{\frac {1}{360}} \left ({b_2}-3+{b_3}-2 {b_1}\right )
\left ({b_4}+3+{b_1}+{b_3}+{b_2}\right )\left ({b_4}+2+
{b_1}+{b_3}+{b_2}\right )$

$\left ({b_4}+1+{b_1}+{b_3}+ {b_2}\right )
(b_2^{2}+2 {b_2} {b_3}+9 {b_2
}+2 {b_1} {b_2}-3 {b_2} {b_4}+9 {b_3}$

$+b_1^{
2}+20+9 {b_1}+2 {b_1} {b_3}+b_3^{2}-3 {b_3} {
b_4}-21 {b_4}-3 {b_1} {b_4}+6 b_4^{2} )$


\item If $ {\rm min} \left \{-{b_2},{b_3}+{b_4},{b_1}+{b_2},-{b_1}-{\it
b_2}-{b_4}\right \} \geq 0  $ then

${\frac {1}{360}} \left ({b_2}+1+{b_1}\right )\left ({b_2}+5+
{b_1}\right )\left ({b_2}+4+{b_1}\right )\left ({b_2}+{
b_1}+3\right )$

$({b_2}+2+{b_1}) (2 {b_2}+
2 {b_1}+3 {b_4}+3+3 {b_3} )$


\item If $ {\rm min} \left \{-{b_1}-{b_4},-{b_2}-{b_3}-{b_4},{b_1}+{\it
b_2}+{b_3}+{b_4},-{b_1}-{b_3}\right \} \geq 0  $ then

${\frac {1}{360}} \left ({b_3}-2 {b_4}+{b_2}+3\right )\left
({b_4}+1+{b_1}+{b_3}+{b_2}\right )\left ({b_4}+3+{\it
b_1}+{b_3}+{b_2}\right )$

$\left ({b_4}+2+{b_1}+{b_3}+{\it
b_2}\right )(b_2^{2}+2 {b_2} {b_4}+2 {b_2} {
b_3}$

$-6 {b_2}-3 {b_1} {b_2}+20+b_4^{2}+2 {b_3
} {b_4}+b_3^{2}-6 {b_4}-3 {b_1} {b_4}-6 {\it
b_3}-3 {b_1} {b_3}+24 {b_1}+6 b_1^{2} )$


\item If $ {\rm min} \left \{{b_1},{b_2},{b_3}+{b_4},-{b_1}-{b_2}-{b_4
}\right \} \geq 0  $ then

${\frac {1}{360}} \left ({b_1}+3\right )\left ({b_1}+2\right )
\left ({b_1}+1\right )\left (b_1^{2}+5 {b_1} {b_2}+9
 {b_1}+20+10 b_2^{2}+30 {b_2}\right )$

$\left (2 {b_2
}+2 {b_1}+3 {b_4}+3+3 {b_3}\right )$


\item If $ {\rm min}\left \{{b_1},{b_2}+{b_3}+{b_4},-{b_1}-{b_3}-{b_4
},-{b_1}-{b_2}-{b_4}\right \} \geq 0  $ then

${\frac {1}{360}} \left ({b_1}+3\right )\left ({b_1}+2\right )
\left ({b_1}+1\right )\left (-{b_4}+3-{b_3}+2 {b_2}
\right )$
$(10 b_2^{2}+30 {b_2}+15 {b_1} {b_2}+
20 {b_2} {b_3}+20 {b_2} {b_4}+20+15 {b_1} {b_4
}+30 {b_3}+30 {b_4}+15 {b_1} {b_3}+6 b_1^{2}+
10 b_4^{2}+10 b_3^{2}+24 {b_1}+20 {b_3} {\it
b_4} )$


\item If $ {\rm min} \left \{-{b_2}-{b_3}-{b_4},-{b_1}-{b_3}-{b_4},{\it
b_1}+{b_2}+{b_3}+{b_4},-{b_1}-{b_2}-{b_4}\right \} \geq 0  $ then

${\frac {1}{360}} \left (-{b_4}+3-{b_3}+2 {b_2}\right )
\left ({b_4}+2+{b_1}+{b_3}+{b_2}\right )\left ({b_4}+1+
{b_1}+{b_3}+{b_2}\right ) $

$\left ({b_4}+3+{b_1}+{b_3}+
{b_2}\right )(b_2^{2}+2 {b_2} {b_4}+2 {b_2
} {b_3}-6 {b_2}-3 {b_1} {b_2}+20+b_4^{2}$

$ +2 {
b_3} {b_4}+b_3^{2}-6 {b_4}-3 {b_1} {b_4}-6
{b_3}-3 {b_1} {b_3}+24 {b_1}+6 b_1^{2} )$

\end{enumerate}
}

\section{Appendix: {\tt Macaulay 2} program}

In this appendix we present our implementation of
the Gr\"obner basis algorithm from Section 2.
For an introduction to the computer algebra system
{\tt Macaulay 2} see  \cite{egss} and \cite{macaulay}.
Our program starts by defining the unimodular
$8 \times 16$-matrix {\tt A} of rank $7$
which represents the counting problem for
$ 4 \times 4$- tables.

\begin{verbatim}
A = {{1, 0, 0, 0,  1, 0, 0, 0,  1, 0, 0, 0,  1, 0, 0, 0},
     {0, 1, 0, 0,  0, 1, 0, 0,  0, 1, 0, 0,  0, 1, 0, 0},
     {0, 0, 1, 0,  0, 0, 1, 0,  0, 0, 1, 0,  0, 0, 1, 0},
     {0, 0, 0, 1,  0, 0, 0, 1,  0, 0, 0, 1,  0, 0, 0, 1},
     {1, 1, 1, 1,  0, 0, 0, 0,  0, 0, 0, 0,  0, 0, 0, 0},
     {0, 0, 0, 0,  1, 1, 1, 1,  0, 0, 0, 0,  0, 0, 0, 0},
     {0, 0, 0, 0,  0, 0, 0, 0,  1, 1, 1, 1,  0, 0, 0, 0},
     {0, 0, 0, 0,  0, 0, 0, 0,  0, 0, 0, 0,  1, 1, 1, 1}}
\end{verbatim}
We next input a  weight vector ${\tt W}$ of length $16$,
to be interpreted as a $4 \times 4$-table:
\begin{verbatim}
W = {1000,1,1,1, 1,1,1,1, 1,1,1,1, 1,1,1,1}
\end{verbatim}

The following five command lines compute the
monomial ideal $\, M = in_w(J_A) $,
here called {\tt nonfaces}, which represents
the chamber we are interested in:
\begin{verbatim}
n = # W;  d = n - # A + 1; R = QQ[x_1..x_n, Weights => W]
Binomial = (b,R) -> (     top := 1_R; bottom := 1_R;
   scan(#b, i -> if b_i > 0 then top = top * R_i^(b_i)
          else bottom = bottom * R_i^(-b_i)); top - bottom);
nonfaces = ideal leadTerm ideal apply( A, a -> Binomial(a,R));
\end{verbatim}
We compute the presentation ideal  $M + L_A$
of the cohomology ring $H^*(X_w;k)$. It is denoted {\tt I}.
\begin{verbatim}
S = QQ[x_1..x_n, r1,r2,r3,r4,c1,c2,c3,c4];
f = map(S,R, toList(x_1..x_n));
Linform = (b,R) -> (s := 0_R; scan(#b,i -> s = s + b_i*R_i); s);
I = f(nonfaces) +
  ideal apply(entries transpose syz matrix A, a -> Linform(a,S));
\end{verbatim}
The next four lines compute a
representation of the Todd class modulo {\tt I}.
\begin{verbatim}
todd = (x) -> (1+1/2*x+1/12*x^2-1/720*x^4+1/30240*x^6-1/1209600*x^8);
trunc = (d,f) -> sum select(terms f, t ->  sum degree t < d+1);
toddclass := 1_S;
scan(1..n, i -> toddclass = trunc(d, toddclass * todd(x_i)) % I);
\end{verbatim}
All subsequent computations take place in the
quotient ring ${\tt T} = {\tt S/I} = $  \break $  H^*(X_w;k)$.
We compute all successive powers of a general divisor $\,\sum u_i x_i $.
\begin{verbatim}
T = S/I;
g = map(T,T,join(toList(n:1) , {r1,r2,r3,r4, c1,c2,c3,c4}));
u = (0, r1-c2-c3-c4,c2,c3,c4,r2,0,0,0,r3,0,0,0,r4,0,0,0);
divp = 1;
divpowers = apply(1..d, i ->
     (divp =  sum toList apply(1..n, i -> u_i * x_i * divp)));
\end{verbatim}
In the final four lines of code, the graded components of
the Todd class are
multiplied with the complementary powers of
the divisor $\,\sum u_i x_i $. The
products are added up (in {\tt T})
and the sum is normalized so that its
constant term is $1$:
\begin{verbatim}
component = (d,f) -> sum select(terms f, t -> d == sum degree t);
erhart = sum toList apply(0..d-1,
  i -> (divpowers_i * (1/(i+1)!) * component(d-i-1,toddclass)));
toString (g(erhart)/g(component(d,toddclass)) + 1)
\end{verbatim}
The final output is a polynomial of degree $9$  in
the variables ${\tt r1}$, ${\tt r2}$, ${\tt r3}$, ${\tt r4}$,
${\tt c1}$, ${\tt c2}$, ${\tt c3}$, ${\tt c4}$. This
particular chamber polynomial has 1967 terms.
The running time of this entire piece of code is about 25 minutes.

Users of {\tt Macaulay 2} will find it easy to modify  our code
so that it works for
any unimodular matrix $A$ and any right hand side $b = Aw $.
Besides redefining the variables {\tt A} and {\tt w}, one only needs
to change those command lines which involve the variables
${\tt r1}$, ${\tt r2}$, ${\tt r3}$, ${\tt r4}$,
${\tt c1}$, ${\tt c2}$, ${\tt c3}$, ${\tt c4}$
particular to $4 \times 4$-tables.

\end{document}